\documentclass[reqno,fleqn]{amsart}%
\usepackage[english]{babel}
\usepackage{amsmath}
\usepackage{amssymb}
\usepackage{amsfonts}
\usepackage{graphicx}%
\theoremstyle{plain}
\newtheorem{definition}{Definition}
\newtheorem{lemma}{Lemma}
\newtheorem{proposition}[lemma]{Proposition}
\newtheorem{theorem}{Theorem}

\theoremstyle{remark}

\newtheorem{example}{Example}

\numberwithin{equation}{section} \allowdisplaybreaks
\begin{document}
\title{The bisymplectomorphism group of a bounded symmetric domain}
\author{Antonio J. Di Scala}
\address{A.D.: Dipartimento di Matematica, Politecnico di Torino, Corso Duca degli
Abruzzi 24, 10129 Torino, Italy }
\email{antonio.discala@polito.it}
\author{Andrea Loi}
\address{A.L.: Dipartimento di Matematica e Informatica, Universit\`{a} di Cagliari,
Via Ospedale 72, 09124 Cagliari, Italy}
\email{loi@unica.it}
\author{Guy Roos}
\address{G.R.: Nevski prospekt 113/4-53, 191024 St Petersburg, Russian Federation}
\email{guy.roos@normalesup.org}
\thanks{Research partially supported by GNSAGA (INdAM) and MIUR of Italy}
\date{February 23, 2008}

\begin{abstract}
We determine the group of diffeomorphisms of a bounded symmetric domain, which
preserve simultaneously the hyperbolic and the flat symplectic form.

\end{abstract}
\maketitle
\tableofcontents

\section*{Introduction}

Let $\Omega$ be an Hermitian bounded symmetric domain in a complex vector
space $V$ ; we always assume that $\Omega$ is given in its circled
realization. The domain $\Omega$ is endowed with two natural symplectic forms:
the flat form $\omega_{0}$ and the hyperbolic form $\omega_{-}$. In a similar
way, the ambient vector space $V$ is also endowed with two natural symplectic
forms: the Fubini-Study form $\omega_{+}$ and the flat form $\omega_{0}$ (see
Section \ref{SEC1} for the definition of $\omega_{-}$, $\omega_{0}$,
$\omega_{+}$). It has been shown in \cite{DiScalaLoi2006} that there exists a
diffeomorphism $F:\Omega\rightarrow V$ such that%
\begin{equation}
F^{\ast}\omega_{0}=\omega_{-},\quad F^{\ast}\omega_{+}=\omega_{0}. \label{I01}%
\end{equation}
This map is the same as the map $\psi$ which was used in \cite[Theorem
VII.4.3]{Roos2000}, with the property%
\[
\psi^{\ast}\omega_{+}^{n}=\omega_{0}^{n}%
\]
($n=\dim_{\mathbb{C}}V$), to show that the flat volume of a bounded symmetric
domain, with some natural normalization, is equal to the degree of a canonical
projective embedding of its compact dual. In the one-dimensional case, where
$V=\mathbb{C}$ and $\Omega$ is the unit disc $\Delta$, this map is simply
$f:\Delta\rightarrow\mathbb{C}$ given by
\begin{equation}
f(z)=\frac{z}{\sqrt{1-\left\vert z\right\vert ^{2}}}. \label{I03}%
\end{equation}
Even in this case, it does not seem that the property (\ref{I01}) had been
noticed before. In the general case, the map $F$ may be defined by
\begin{equation}
F(z)=B(z,z)^{-1/4}z, \label{I04}%
\end{equation}
where $B(z,z)$ denotes the Bergman operator of the Jordan triple structure on
$V$ associated to $\Omega$; it may also be defined by functional calculus in
Hermitian positive Jordan triples. In view of the property (\ref{I01}), the
map $F$ is called map of (bi)symplectic duality. The part $F^{\ast}\omega
_{0}=\omega_{-}$ tells that $F$ is a realization of the isomorphism of Mc Duff
\cite{McDuff1988} for the bounded symmetric domain $\Omega$; but the property
(\ref{I01}), which involves two pairs of symplectic forms, is much stronger.
In order to determine all diffeomorphisms $F:\Omega\rightarrow V$ verifying
(\ref{I01}), we determine the group of bisymplectomorphims of $\Omega$, that
is, diffeomorphisms $\phi:\Omega\rightarrow\Omega$ such that%
\begin{equation}
\phi^{\ast}\omega_{0}=\omega_{0},\quad\phi^{\ast}\omega_{-}=\omega_{-}.
\label{I02}%
\end{equation}
This group is infinite-dimensional, but is the direct product of the
compact group $K$ of linear automorphisms of $\Omega$ with an
infinite-dimensional Abelian group of \textquotedblleft radial circular
diffeomorphisms\textquotedblright\ (Theorem \ref{structure}); this is the main
result of this article and may be considered as a kind of Schwarz lemma.

The plan of this article is as follows. In Section \ref{SEC1}, we recall known
facts about Jordan triple systems associated to bounded complex symmetric
domains (see mainly \cite{Loos1977}); the only result we could not find in the
literature is Proposition \ref{frames}, which describes the tangent space of
the manifold of frames (the \textquotedblleft F\"{u}rstenberg-Satake
boundary\textquotedblright\ of $\Omega$) in terms of Peirce decomposition in
Jordan triples. In Section \ref{SEC2}, we compute the symplectic forms
$\omega_{0}$ and $\omega_{-}$ using spectral decomposition in Jordan triples,
which is the appropriate generalization of polar coordinates. From this, we
derive in Section \ref{SEC3} a simple proof of the property (\ref{I01}),
different from the proof given in \cite{DiScalaLoi2006}. Section \ref{SEC4} is
devoted to the study and characterization of bisymplectomorphisms.

%\bigskip
%The authors want to thank the referee for useful remarks.

\section{Hermitian positive Jordan triples\label{SEC1}}

Let $\Omega$ be a bounded symmetric domain in a finite dimensional complex
vector space $V$. \emph{We will always consider such a domain in its (unique
up to linear isomorphism) circled realization.} Consider the associated Jordan
triple system $\left(  V,\left\{  ~,~,~\right\}  \right)  $. For basic facts
about Hermitian positive Jordan triples and their correspondence with complex
symmetric domains, see \cite{Loos1977}, \cite{Roos2000}. We recall hereunder
those which will be used here.

\subsection{Definitions and notations}

Consider the operators on the Jordan triple $V$ defined by
\begin{align}
D(x,y)z  &  =\left\{  x,y,z\right\}  ,\label{D1}\\
Q(x,z)y  &  =\left\{  x,y,z\right\}  ,\label{D2}\\
Q(x,x)  &  =2Q(x),\label{D3}\\
B(x,y)  &  =\operatorname{id}_{V}-D(x,y)+Q(x)Q(y). \label{D4}%
\end{align}
The operators $D(x,y)$ and $B(x,y)$ are $\mathbb{C}$-linear, the operator
$Q(x)$ is $\mathbb{C}$-antilinear. The hermitian form%
\begin{equation}
\left(  u\mid v\right)  =\operatorname{tr}D(u,v) \label{D5}%
\end{equation}
is a Hermitian scalar product on $V$; with respect to this product, $D(x,x)$
and $B(x,x)$ are self-adjoint.

For $z\in V$, the \emph{odd powers} $z^{(2p+1)}$ of $z$ in the Jordan triple
system $V$ are defined by
\begin{equation}
z^{(1)}=z,\qquad z^{(2p+1)}=Q(z)z^{(2p-1)}. \label{D6}%
\end{equation}
An element $e\in V$ is called \emph{tripotent} if $e\neq0$ and $e^{(3)}=e$.
Two tripotents $e_{1},e_{2}$ are called \emph{(strongly) orthogonal} if
$D\left(  e_{1},e_{2}\right)  =0$. A tripotent element is called
\emph{minimal}, or \emph{primitive}, if it is not the sum of two orthogonal
tripotents. A tripotent element $e$ is called \emph{maximal} if there is no
tripotent orthogonal to $e$.

\subsection{Spectral decomposition}

Each element $z\in V$ has a unique \emph{spectral decomposition}%
\begin{equation}
z=\lambda_{1}e_{1}+\cdots+\lambda_{s}e_{s}\qquad(\lambda_{1}>\cdots
>\lambda_{s}>0), \label{D7}%
\end{equation}
where $\left(  e_{1},\ldots,e_{s}\right)  $ is a sequence of pairwise
orthogonal tripotents. The integer $s=\operatorname{rk}z$ is called the
\emph{rank} of $z$. Let $V_{z}^{+}$ be the $\mathbb{R}$-subspace of $V$
generated by the odd powers $z,\ldots,z^{(2p+1)},\ldots$ and $V_{z}=V_{z}%
^{+}\oplus\operatorname{i}V_{z}^{+}$ the $\mathbb{C}$-subspace generated by
the odd powers of $z$. Then $\operatorname{rk}z=\dim_{\mathbb{R}}V_{z}^{+}$
and $\left(  e_{1},\ldots,e_{s}\right)  $ is an $\mathbb{R}$-basis of
$V_{z}^{+}$. The \emph{rank} of $V$ is $r=\operatorname{rk}V=\max\left\{
\operatorname{rk}z\mid z\in V\right\}  $; elements $z$ such that
$\operatorname{rk}z=\operatorname{rk}V$ are called \emph{regular}. If $z\in V$
is regular, with spectral decomposition%
\begin{equation}
z=\lambda_{1}e_{1}+\cdots+\lambda_{r}e_{r}\qquad(\lambda_{1}>\cdots
>\lambda_{r}>0), \label{D8}%
\end{equation}
then $\left(  e_{1},\ldots,e_{r}\right)  $ is a \emph{(Jordan) frame of }$V$,
that is, a maximal sequence of pairwise orthogonal minimal tripotents.

\subsection{Peirce decompositions}

Let $e\in V$ be a tripotent. Then the eigenvalues of $D(e,e)$ are contained in
$\left\{  0,1,2\right\}  $. Define the \emph{Peirce subspaces} of $e$ as
\begin{equation}
V_{i}(e)=\left\{  z\in V\mid D(e,e)z=iz\right\}  \qquad(i\in\left\{
0,1,2\right\}  ). \label{P1}%
\end{equation}
The decomposition
\begin{equation}
V=V_{0}(e)\oplus V_{1}(e)\oplus V_{2}(e) \label{P2}%
\end{equation}
is called the \emph{Peirce decomposition} of $V$ w.r.~to $e$. A tripotent $e$
is maximal if $V_{0}(e)=0$, minimal if $V_{2}(e)=\mathbb{C}e$. The Peirce
subspaces compose according to the law%
\begin{equation}
\left\{  V_{i}(e),V_{j}(e),V_{k}(e)\right\}  \subset V_{i-j+k}(e), \label{P3}%
\end{equation}
where $V_{m}(e)=0$ if $m\notin\left\{  0,1,2\right\}  $; in particular, Peirce
subspaces are Jordan subsystems of $V$. The $\mathbb{C}$-antilinear operator
$Q(e)$ is $0$ on $V_{0}(e)\oplus V_{1}(e)$; its restriction to $V_{2}(e)$ is
involutive. Let
\begin{equation}
V_{2}^{+}(e)=\left\{  v\in V\mid D(e,e)v=2v,\ Q(e)v=v\right\}  . \label{P5}%
\end{equation}
Then the decomposition of $V(e)$ into (real) eigenspaces of $Q(e)$ is%
\begin{equation}
V_{2}(e)=V_{2}^{+}(e)\oplus\operatorname{i}V_{2}^{+}(e). \label{P6}%
\end{equation}

Let $\mathbf{e}=\left(  e_{1},\ldots,e_{s}\right)  $ be a sequence of pairwise
orthogonal tripotents. Then the operators $D(e_{j},e_{j})$, $1\leq j\leq s$
commute and have the common eigenspaces%
\begin{align}
V_{jj}(\mathbf{e})  &  =V_{2}(e_{j})\qquad(1\leq j\leq s),\nonumber\\
V_{jk}(\mathbf{e)}  &  =V_{1}(e_{j})\cap V_{1}(e_{k})\qquad(1\leq j<k\leq
s),\nonumber\\
V_{0j}(\mathbf{e})  &  =V_{1}(e_{j})\cap{\displaystyle\bigcap\limits_{k\neq
j}} V_{0}(e_{k})\qquad(1\leq j\leq s),\label{P7}\\
V_{00}(\mathbf{e})  &  ={\displaystyle\bigcap\limits_{k}} V_{0}(e_{k}%
)\nonumber
\end{align}
(some of these subspaces may be $0$). The decomposition
\begin{equation}
V={\displaystyle\bigoplus\limits_{0\leq j\leq k\leq s}} V_{jk}(\mathbf{e})
\label{P8}%
\end{equation}
is called the \emph{simultaneous Peirce decomposition} of $V$ w.r.\ to
$\mathbf{e}$. If
\begin{align*}
z  &  =\lambda_{1}e_{1}+\cdots+\lambda_{s}e_{s}\qquad(\lambda_{j}\in
\mathbb{C}),\\
e  &  =e_{1}+\cdots+e_{s}%
\end{align*}
and $v\in V_{jk}(\mathbf{e})$, then
\begin{align}
D(z,z)v  &  =\left(  \left\vert \lambda_{j}\right\vert ^{2}+\left\vert
\lambda_{k}\right\vert ^{2}\right)  v\label{P9}\\
Q(z)v  &  =\lambda_{j}\lambda_{k}Q(e)v,\label{P10}\\
Q(z)Q(z)  &  =\left\vert \lambda_{j}\lambda_{k}\right\vert ^{2}v,\label{P11}\\
B(z,z)v  &  =\left(  1-\left\vert \lambda_{j}\right\vert ^{2}\right)  \left(
1-\left\vert \lambda_{k}\right\vert ^{2}\right)  v,\label{P12}\\
B(z,-z)v  &  =\left(  1+\left\vert \lambda_{j}\right\vert ^{2}\right)  \left(
1+\left\vert \lambda_{k}\right\vert ^{2}\right)  v, \label{P12+}%
\end{align}
where $\lambda_{0}=0$. So the $V_{jk}(\mathbf{e})$'s are eigenspaces for all
the operators $D(z,z)$, $B(z,z)$, $B(z,-z)$, $z=\lambda_{1}e_{1}%
+\cdots+\lambda_{s}e_{s}$.

\subsection{Hermitian metrics and symplectic forms}

Let $V$ be a Hermitian positive Jordan triple and let $\Omega$ be the
associated Hermitian bounded symmetric domain. Let
\begin{equation}
h_{0}(z)(u,v)=\left(  u\mid v\right)  =\operatorname{tr}D(u,v) \label{G02}%
\end{equation}
be the flat Hermitian metric and let
\begin{align*}
\omega_{0}(z)  &  =\textstyle\frac{\operatorname{i}}{2}\partial\overline{\partial
}\left(  z\mid z\right)  ,\\
\omega_{0}(z)(u,v)  &  =\textstyle\frac{\operatorname{i}}{2}\left(  \left(  u\mid
v\right)  -\left(  v\mid u\right)  \right)
\end{align*}
be the associated\emph{ flat symplectic form}. If $\Omega$ is endowed with the
volume form $\omega_{0}^{n}$ ($n=\dim_{\mathbb{C}}V$), the Bergman kernel of
$\Omega$ is
\[
K(x,y)=\frac{C}{\det B(x,y)},
\]
with $C=\left(  \int_{\Omega}\omega_{0}^{n}\right)  ^{-1}$. The \emph{Bergman
metric} of $\Omega$ is
\[
h_{-}(z)(u,v)=\partial_{u}\overline{\partial}_{v}\ln K(z,z)=-\partial
_{u}\overline{\partial}_{v}\ln\det B(z,z).
\]
It satisfies the relation
\begin{equation}
h_{-}(z)(u,v)=h_{0}\left(  B(z,z)^{-1}u,v\right)  . \label{G01}%
\end{equation}
In view of this relation, $B(z,z)$ is called the \emph{Bergman operator} at
$z\in\Omega$.

The \emph{hyperbolic symplectic form} of $\Omega$, associated to the Bergman
metric, is defined by
\begin{equation}
\omega_{-}(z)=-\textstyle\frac{\operatorname{i}}{2}\partial\overline{\partial}\ln\det
B(z,z). \label{G03}%
\end{equation}
From (\ref{G01}), it results that the forms $\omega_{0}$ and $\omega_{-}$ are
related with the Bergman operator by%
\begin{equation}
\omega_{-}(z)(u,v)=\omega_{0}(B(z,z)^{-1}u,v), \label{G05}%
\end{equation}
for $z\in\Omega$ and $u,v\in T_{z}\Omega$.

The (generalized) \emph{Fubini-Study} metric on $V$ is defined by%
\[
h_{+}(z)(u,v)=\partial_{u}\overline{\partial}_{v}\ln\det B(z,-z).
\]
The associated K\"{a}hler form is
\[
\omega_{+}(z)=\textstyle\frac{\operatorname{i}}{2}\partial\overline{\partial}\ln\det
B(z,-z).
\]
It is related to the flat form by%
\begin{equation}
\omega_{+}(z)(u,v)=\omega_{0}(B(z,-z)^{-1}u,v). \label{G04}%
\end{equation}

\subsection{Polar coordinates}

Let $M$ be the set ot tripotents elements of the positive Jordan triple $V$.
Then $M$ is a compact submanifold of $V$ (with connected components of
different dimensions). At $e\in M$, the tangent space $T_{e}M$ and the normal
space $N_{e}M$ to $M$ are%
\begin{align}
T_{e}M  &  =\operatorname{i}V_{2}^{+}(e)\oplus V_{1}(e),\label{P13}\\
N_{e}M  &  =V_{0}(e)\oplus V_{2}^{+}(e) \label{P14}%
\end{align}
(see \cite[Theorem 5.6]{Loos1977}).

The \emph{height} $k$ of a tripotent element $e$ is the maximal length of a
decomposition $e=e_{1}+\cdots+e_{k}$ into a sum of pairwise orthogonal
(minimal) tripotents. Minimal tripotents have height $1$, maximal tripotents
have height $r=\operatorname{rk}V$. Denote by $M_{k}$ the set of tripotents of
height $k$. \emph{If }$V$ \emph{is simple (that is, if }$\Omega$ \emph{is
irreducible)}, the submanifolds $M_{k}$ are the connected components of $M$.

The set $\mathcal{F}$ of frames (also called F\"{u}rstenberg-Satake boundary
of $\Omega$):
\begin{equation}
\mathcal{F}=\left\{  \left(  e_{1},\ldots,e_{r}\right)  \mid e_{j}\in
M_{1},\ e_{j}\perp e_{k}\ (1\leq j<k\leq r)\right\}  ,\label{P15}%
\end{equation}
(where $e_{j}\perp e_{k}$ means orthogonality of tripotents: $D(e_{j}%
,e_{k})=0$ or equivalently $\left\{  e_{j},e_{j},e_{k}\right\}  =0$) is a
submanifold of $V^{r}$. If $\mathbf{e=}\left(  e_{1},\ldots,e_{r}\right)  $ is
a frame for a simple positive Hermitian Jordan triple $V$, the corresponding
Peirce subspaces have the following properties (see \cite[Theorem
VI.3.5]{Roos2000}):

\begin{itemize}
\item $V_{00}(\mathbf{e})=0$;

\item $V_{jj}(\mathbf{e})=\mathbb{C}e_{j}$ ($1\leq j\leq r$);

\item all $V_{jk}(\mathbf{e})=V_{kj}(\mathbf{e})$ ($1\leq j<k\leq r$) have the
same dimension $a>0$;

\item all $V_{0j}(\mathbf{e})$ ($1\leq j\leq r$) have the same dimension $b$;
these subspaces are $0$ if and only if the domain $\Omega$ is of
\textquotedblleft tube type\textquotedblright.
\end{itemize}

The following proposition provides a description of the tangent space of
$\mathcal{F}$.

\begin{proposition}
\label{frames}Let $\mathbf{e=}\left(  e_{1},\ldots,e_{r}\right)
\in\mathcal{F}\subset V^{r}$ and $e=e_{1}+\cdots+e_{r}$. Then $\left(
v_{1},\ldots,v_{r}\right)  \in T_{\mathbf{e}}\mathcal{F}$ if and only if
\begin{equation}
v_{j}=\operatorname{i}\alpha_{j}e_{j}+v_{j0}+\sum_{\substack{1\leq k\leq
r\\k\neq j}}v_{jk}\qquad(1\leq j\leq r), \label{P17}%
\end{equation}
where $\alpha_{j}\in\mathbb{R}$, $v_{j0}\in V_{0j}(\mathbf{e})$, $v_{jk}\in
V_{jk}(\mathbf{e})=V_{kj}(\mathbf{e})$ and
\begin{equation}
Q(e)v_{jk}=-v_{kj}\qquad(1\leq j<k\leq r). \label{P19}%
\end{equation}

\end{proposition}

\begin{proof}
Let $\left(  v_{1},\ldots,v_{r}\right)  \in T_{\mathbf{e}}\mathcal{F}$. As
$e_{j}$ are minimal tripotents, we have%
\[
v_{j}\in T_{e_{j}}M_{1}=\operatorname{i}\mathbb{R}e_{j}\oplus V_{1}%
(e_{j})=\operatorname{i}\mathbb{R}e_{j}\oplus V_{0j}(\mathbf{e})\oplus
{\bigoplus\limits_{\substack{1\leq k\leq r\\k\neq j}}}V_{jk}(\mathbf{e}),
\]
which shows that $v_{j}$ has the form (\ref{P17}).

The orthogonality conditions in a frame are
\[
\left\{  e_{j},e_{j},e_{k}\right\}  =0\qquad(1\leq j<k\leq r).
\]
Differentiating these conditions yields%
\[
\left\{  v_{j},e_{j},e_{k}\right\}  +\left\{  e_{j},v_{j},e_{k}\right\}
+\left\{  e_{j},e_{j},v_{k}\right\}  =0\qquad(1\leq j<k\leq r).
\]
As $D(e_{j},e_{k})=0$, this condition is reduced to
\begin{equation}
Q(e_{j},e_{k})v_{j}+D(e_{j},e_{j})v_{k}=0\qquad(1\leq j<k\leq r). \label{P18}%
\end{equation}
Let
\[
v_{j}=\operatorname{i}\alpha_{j}e_{j}+v_{j0}+\sum_{\substack{1\leq m\leq
r\\m\neq j}}v_{jm}\qquad(1\leq j\leq r).
\]
Then%
\begin{align*}
D(e_{j},e_{j})v_{k}  &  =v_{kj},\quad Q(e_{j},e_{k})e_{j}=0,\\
Q(e_{j},e_{k})  &  =Q(e_{j}+e_{k})-Q(e_{j})-Q(e_{k}),
\end{align*}
and we get from (\ref{P10})%
\begin{align*}
Q(e_{j}+e_{k})v_{jm}  &  =\delta_{k}^{m}Q(e)v_{jm},\quad Q(e_{j}%
)v_{jm}=0,\quad Q(e_{k})v_{jm}=0,\\
Q(e_{j},e_{k})v_{j}  &  =Q(e)v_{jk}.
\end{align*}
This shows that the conditions (\ref{P18}) are equivalent to (\ref{P19}). 
\end{proof}

Comparing the description of $T_{\mathbf{e}}\mathcal{F}$ with the simultaneous
Peirce decomposition of $V$ w.r.~to $\mathbf{e}$, it is easily checked that
$T_{\mathbf{e}}\mathcal{F}$ is a real vector space of dimension $2n-r$, where
$n=\dim_{\mathbb{C}}V$. This implies that the map
\begin{align}
\left\{  \lambda_{1}>\cdots>\lambda_{r}>0\right\}  \times\mathcal{F}  &
\rightarrow V_{\mathrm{reg}}\nonumber\\
\left(  \left(  \lambda_{1},\ldots,\lambda_{r}\right)  ,\left(  e_{1}%
,\ldots,e_{r}\right)  \right)   &  \mapsto\sum\lambda_{j}e_{j} \label{P16}%
\end{align}
is a diffeomorphism onto the set $V_{\mathrm{reg}}$ of regular elements of
$V$; its restriction%
\[
\left\{  1>\lambda_{1}>\cdots>\lambda_{r}>0\right\}  \times\mathcal{F}%
\rightarrow\Omega_{\mathrm{reg}}%
\]
is a diffeomorphism onto the set $\Omega_{\mathrm{reg}}$ of regular elements
of $\Omega$. This map plays the same role as polar coordinates in rank one.

\subsection{Functional calculus\label{FC}}

See \cite[Section 3.18]{Loos1977}. Using the spectral decomposition, it is
possible to associate to an \emph{odd} function $f:(-1,1)\rightarrow
\mathbb{C}$ (resp. $f:\mathbb{R}\rightarrow\mathbb{C}$) a \textquotedblleft
radial\textquotedblright\ map $F:\Omega\rightarrow V$ (resp. $F:V\rightarrow
V$) in the following way. Let $z\in V$ and let%
\[
z=\lambda_{1}e_{1}+\cdots+\lambda_{k}e_{k},\quad\lambda_{1}>\cdots>\lambda
_{k}>0
\]
be the spectral decomposition of $z$. Define the map $F=\widetilde{f}$
associated to $f$ by
\begin{equation}
F(z)=f(\lambda_{1})e_{1}+\cdots+f(\lambda_{k})e_{k}. \label{F00}%
\end{equation}
Since $f$ is odd, it follows from properties of tripotents that, for
\[
z=\lambda_{1}e_{1}+\cdots+\lambda_{r}e_{r},
\]
where $\mathbf{e}=(e_{1},\ldots,e_{r})$ is a frame and $\lambda_{1}%
,\ldots,\lambda_{r}\in\mathbb{R}$, we have
\begin{equation}
F(z)=f(\lambda_{1})e_{1}+\cdots+f(\lambda_{r})e_{r}. \label{F01}%
\end{equation}
The linear correspondence $f\mapsto F$ has the following properties:

\begin{enumerate}
\item If $f$ is continuous, then $F$ is continuous.

\item If
\[
f(t)=\sum_{k=0}^{N}a_{k}t^{2k+1}%
\]
is a polynomial, then $F$ is the map defined by%
\begin{equation}
F(z)=\sum_{k=0}^{N}a_{k}z^{(2k+1)}\qquad(z\in V). \label{F02}%
\end{equation}

\item If $f$ is analytic, then $F$ is real-analytic; if $f$ is given near $0$
by%
\[
f(t)=\sum_{k=0}^{\infty}a_{k}t^{2k+1},
\]
then $F$ has the Taylor expansion near $0\in V$:%
\begin{equation}
F(z)=\sum_{k=0}^{\infty}a_{k}z^{(2k+1)}.\label{F03}%
\end{equation}

\item If $\left\vert f(\lambda)\right\vert \leq C\left\vert \lambda\right\vert
^{2N+1}$ for $\left\vert \lambda\right\vert <c$, then $\left\Vert
F(z)\right\Vert \leq C\left\Vert z\right\Vert ^{2N+1}$ for $z\in V$,
$\left\Vert z\right\Vert <c$. (Here $\left\Vert z\right\Vert $ denotes the
\emph{spectral norm} of $z$ in $V$: $\left\Vert z\right\Vert =\max\left\vert
\lambda_{j}\right\vert $ for $z=\lambda_{1}e_{1}+\cdots+\lambda_{r}e_{r}$,
where $\mathbf{e}=(e_{1},\ldots,e_{r})$ is a frame).

\item If $f$ is $C^{\infty}$, then $F$ is also $C^{\infty}$.
\end{enumerate}

\noindent Properties (1)--(4) are easy. Property (5) follows from (2) and (4).

\section{Symplectic forms in polar coordinates\label{SEC2}}

\subsection{The flat symplectic form}

We compute the flat symplectic form $\omega_{0}$ in the \textquotedblleft
polar coordinates\textquotedblright\
\begin{align*}
\left\{  \lambda_{1}>\cdots>\lambda_{r}>0\right\}  \times\mathcal{F} &
\rightarrow V_{\mathrm{reg}}\\
\left(  \left(  \lambda_{1},\ldots,\lambda_{r}\right)  ,\left(  e_{1}%
,\ldots,e_{r}\right)  \right)   &  \mapsto\sum\lambda_{j}e_{j}.
\end{align*}
If $\left(  z^{1},\ldots,z^{n}\right)  $ ($n=\dim V$) are orthonormal
coordinates for the Hermitian product $\left(  u\mid v\right)  $, then%
\[
\omega_{0}={\textstyle\frac{\operatorname{i}}{2}}\sum_{m=1}^{n}\operatorname{d}%
z^{m}\wedge\operatorname{d}\overline{z}^{m}.
\]
Let $\left(  e_{j}^{m}\right)  _{1\leq m\leq n}$ be the coordinates of $e_{j}$
in the same basis as for $z$. From%
\[
z=\sum_{j=1}^{r}\lambda_{j}e_{j},\quad z^{m}=\sum_{j=1}^{r}\lambda_{j}%
e_{j}^{m},
\]
we have
\begin{align*}
\operatorname{d}z^{m}\wedge\operatorname{d}\overline{z}^{m} &  =\sum
_{j,k=1}^{r}\lambda_{j}\lambda_{k}\operatorname{d}e_{j}^{m}\wedge
\operatorname{d}\overline{e}_{k}^{m}+\sum_{j,k=1}^{r}e_{j}^{m}\overline{e}%
_{k}^{m}\operatorname{d}\lambda_{j}\wedge\operatorname{d}\lambda_{k}\\
&  \quad+\sum_{j,k=1}^{r}\lambda_{k}\operatorname{d}\lambda_{j}\wedge\left(
e_{j}^{m}\operatorname{d}\overline{e}_{k}^{m}-\overline{e}_{k}^{m}%
\operatorname{d}e_{j}^{m}\right)  .
\end{align*}
Using $\left(  e_{j}\mid e_{k}\right)  =\delta_{jk}$, we get
\begin{equation}
\omega_{0}=\sum_{j,k=1}^{r}\lambda_{j}\lambda_{k}\omega_{jk}+ \sum_{j,k=1}%
^{r}\lambda_{k}\operatorname{d}\lambda_{j}\wedge\eta_{jk},\label{S01}%
\end{equation}
where%
\begin{align}
\eta_{jk} & ={\textstyle\frac{\operatorname{i}}{2}}\left.  \sum_{m=1}^{n}e_{j}%
^{m}\operatorname{d}\overline{e}_{k}^{m}\right\vert _{\mathcal{F}}%
,\label{S07}\\
\omega_{jk} &  ={\textstyle\frac{\operatorname{i}}{2}}\left.  \sum_{m=1}%
^{n}\operatorname{d}e_{j}^{m}\wedge\operatorname{d}\overline{e}_{k}%
^{m}\right\vert _{\mathcal{F}}=\operatorname{d}\eta_{jk}.\label{S08}%
\end{align}

We compute $\eta_{jk}$ and $\omega_{jk}$ using the description of
$T_{\mathbf{e}}\mathcal{F}$ given in Proposition \ref{frames}. Let $v,w\in
T_{\mathbf{e}}\mathcal{F}$ with $v=\left(  v_{1},\ldots,v_{r}\right)  $,
$w=\left(  w_{1},\ldots,w_{r}\right)  $,
\begin{align*}
v_{j} &  =\operatorname{i}\alpha_{j}e_{j}+v_{j0}+\sum_{\substack{1\leq m\leq
r\\m\neq j}}v_{jm}\qquad(1\leq j\leq r),\\
w_{j} &  =\operatorname{i}\beta_{j}e_{j}+w_{j0}+\sum_{\substack{1\leq m\leq
r\\m\neq j}}w_{jm}\qquad(1\leq j\leq r),
\end{align*}
where
\begin{align*}
\alpha_{j},\beta_{j}  & \in\mathbb{R},\quad v_{j0},w_{j0}\in V_{0j}%
(\mathbf{e}),\quad v_{jk},w_{jk}\in V_{jk}(\mathbf{e}),\\
Q(e)v_{jk}  & =-v_{kj},\quad Q(e)w_{jk}=-w_{kj}\qquad(1\leq j<k\leq r).
\end{align*}
Then, as the Peirce subspaces are orthogonal w.r.~to $\left(  \ \mid\ \right)
$, we deduce from (\ref{S07})-(\ref{S08})
\begin{align}
\eta_{jj}(\mathbf{e})(v)  &  ={\textstyle\frac{\operatorname{i}}{2}}\left(  e_{j}\mid v_{j}\right)
={\textstyle\frac{{1}}{2}}\alpha_{j},\label{S06}\\
\eta_{jk}(\mathbf{e})(v)  &  ={\textstyle\frac{\operatorname{i}}{2}}\left(  e_{j}\mid v_{k}\right)
=0\qquad(j\neq k),\nonumber\\
\omega_{jk}  &  =\operatorname{d}\eta_{jk}=0\qquad(j\neq k)\nonumber
\end{align}
and
\begin{align}
\omega_{jj}(\mathbf{e})(v,w)  &  =\textstyle\frac{\operatorname{i}}{2}\left(  \left(
v_{j}\mid w_{j}\right)  -\left(  w_{j}\mid v_{j}\right)  \right) \nonumber\\
&  =\left\langle v_{j0}\mid w_{j0}\right\rangle +\sum_{\substack{1\leq m\leq
r\\m\neq j}}\left\langle v_{jm}\mid w_{jm}\right\rangle , \label{S05}%
\end{align}
where $\left\langle \ \mid\ \right\rangle $ denotes the symplectic product
\begin{equation}
\left\langle x\mid y\right\rangle =\textstyle\frac{\operatorname{i}}{2}\left(  \left(
x\mid y\right)  -\left(  y\mid x\right)  \right)  . \label{S02}%
\end{equation}
Finally, we have
\begin{equation}
\omega_{0}=\sum_{j=1}^{r}\lambda_{j}^{2}\omega_{jj}+ \sum_{j=1}^{r}\lambda
_{j}\operatorname{d}\lambda_{j}\wedge\eta_{jj}, \label{S04}%
\end{equation}
with $\omega_{jj}$ and $\eta_{jj}$ given by (\ref{S05}), (\ref{S06}). The
expression (\ref{S05}) shows that the $\omega_{jj}$'s ($1\leq j\leq r$) are
linearly independent at each point $\mathbf{e}\in\mathcal{F}$.

As
\begin{align*}
2\left(  Q(z)x\mid y\right)   &  =\left(  D(z,x)z\mid y\right)  =\left(  z\mid
D(x,z)y\right) \\
&  =\left(  z\mid D(y,z)x\right)  =\left(  D(z,y)z\mid x\right)  =2\left(
Q(z)y\mid x\right)  ,
\end{align*}
the $Q$ operator satisfies
\begin{align}
\left(  Q(z)x\mid y\right)   &  =\left(  Q(z)y\mid x\right)  ,\label{S03}\\
\left\langle Q(z)x\mid y\right\rangle  &  =-\left\langle x\mid
Q(z)y\right\rangle \label{S09}%
\end{align}
for all $z,x,y\in V$.

For $v,w\in T_{\mathbf{e}}\mathcal{F}$, using
\[
Q(e)v_{jk}=-v_{kj},\quad Q(e)w_{jk}=-w_{kj}\qquad(1\leq j\neq k\leq r),
\]
we have then
\[
\left\langle v_{jk}\mid w_{jk}\right\rangle =\left\langle Q(e)v_{kj}\mid
Q(e)w_{kj}\right\rangle =-\left\langle v_{kj}\mid Q(e)^{2}w_{kj}\right\rangle
,
\]
that is,
\begin{equation}
\left\langle v_{jk}\mid w_{jk}\right\rangle =-\left\langle v_{kj}\mid
w_{kj}\right\rangle , \label{S10}%
\end{equation}
as $w_{kj}\in V_{2}(e)$ and $Q(e)$ is involutive on $V_{2}(e)$. In view of
(\ref{S10}), the flat symplectic form $\omega_{0}$ in polar coordinates may be
rewritten%
\begin{equation}
\omega_{0}=\sum_{j=1}^{r}\lambda_{j}^{2}\theta_{j0}+\sum_{\substack{j,k\\1\leq
j<k\leq r}}\left(  \lambda_{j}^{2}-\lambda_{k}^{2}\right)  \theta_{jk}%
+ \sum_{j=1}^{r}\lambda_{j}\operatorname{d}\lambda_{j}\wedge\eta_{jj},
\label{S11}%
\end{equation}
where the $\eta_{jj}$'s are defined by (\ref{S06}), and the $\theta_{j0}$'s
and $\theta_{jk}$'s by
\begin{align}
\theta_{j0}(\mathbf{e})(v,w)  &  =\left\langle v_{j0}\mid w_{j0}\right\rangle
,\qquad(1\leq j\leq r),\label{S12}\\
\theta_{jk}(\mathbf{e})(v,w)  &  =\left\langle v_{jk}\mid w_{jk}\right\rangle
\qquad(1\leq j<k\leq r). \label{S13}%
\end{align}
for $v,w\in T_{\mathbf{e}}\mathcal{F}$. Note that these forms are (pull-backs
of) forms on the manifold of frames $\mathcal{F}$, and that the Peirce subspaces
$V_{0j}(\mathbf{e})$, hence the $\theta_{j0}$'s, are $0$ when the domain is of tube type.

\subsection{Hyperbolic form and Fubini--Study form}

We compute now the hyperbolic form $\omega_{-}$ and the Fubini--Study form
$\omega_{+}$ in polar coordinates.

Using
\begin{align*}
B(z,z)^{-1}v_{j0}  &  =\left(  1-\lambda_{j}^{2}\right)  ^{-1}v_{j0}%
\qquad(1\leq j\leq r),\\
B(z,z)^{-1}v_{jk}  &  =\left(  1-\lambda_{j}^{2}\right)  ^{-1}\left(
1-\lambda_{k}^{2}\right)  ^{-1}v_{jk}\qquad(1\leq j<k\leq r),\\
B(z,z)^{-1}e_{j}  &  =\left(  1-\lambda_{j}^{2}\right)  ^{-2}e_{j}\qquad(1\leq
j\leq r)
\end{align*}
in (\ref{S11}), (\ref{S12}), (\ref{S13}), we obtain
\begin{align}
\eta_{jj}(\mathbf{e})(B(z,z)^{-1}v)  &  =\left(  1-\lambda_{j}^{2}\right)
^{-2}\eta_{jj}(\mathbf{e})(v),\label{S20}\\
\theta_{j0}(\mathbf{e})(B(z,z)^{-1}v,w)  &  =\left(  1-\lambda_{j}^{2}\right)
^{-1}\theta_{j0}(\mathbf{e})(v,w),\label{S21}\\
\theta_{jk}(\mathbf{e})(B(z,z)^{-1}v,w)  &  =\left(  1-\lambda_{j}^{2}\right)
^{-1}\left(  1-\lambda_{k}^{2}\right)  ^{-1}\theta_{jk}(\mathbf{e})(v,w),
\label{S22}%
\end{align}
for $z=\lambda_{1}e_{1}+\cdots+\lambda_{r}e_{r}$ and $v,w\in T_{\mathbf{e}%
}\mathcal{F}$. From (\ref{G05}):%
\[
\omega_{-}(z)(v,w)=\omega_{0}(B(z,z)^{-1}v,w)
\]
and the expression of $\omega_{0}$ in polar coordinates, we have%
\begin{equation}
\omega_{-}=\sum_{j=1}^{r}\frac{\lambda_{j}^{2}}{1-\lambda_{j}^{2}}\theta
_{j0}+\sum_{\substack{j,k\\1\leq j<k\leq r}}\frac{\lambda_{j}^{2}-\lambda
_{k}^{2}}{\left(  1-\lambda_{j}^{2}\right)  \left(  1-\lambda_{k}^{2}\right)
}\theta_{jk}+ \sum_{j=1}^{r}\frac{\lambda_{j}\operatorname{d}\lambda_{j}%
}{\left(  1-\lambda_{j}^{2}\right)  ^{2}}\wedge\eta_{jj}. \label{S14}%
\end{equation}

In the same way, the Fubini-Study symplectic form on $V$ is%
\begin{equation}
\omega_{+}=\sum_{j=1}^{r}\frac{\lambda_{j}^{2}}{1+\lambda_{j}^{2}}\theta
_{j0}+\sum_{\substack{j,k\\1\leq j<k\leq r}}\frac{\lambda_{j}^{2}+\lambda
_{k}^{2}}{\left(  1+\lambda_{j}^{2}\right)  \left(  1+\lambda_{k}^{2}\right)
}\theta_{jk}+ \sum_{j=1}^{r}\frac{\lambda_{j}\operatorname{d}\lambda_{j}%
}{\left(  1+\lambda_{j}^{2}\right)  ^{2}}\wedge\eta_{jj}. \label{S17}%
\end{equation}

\section{Symplectic duality\label{SEC3}}

Consider the real analytic maps $f=]-1,1[\rightarrow\mathbb{R}$ and
$g:\mathbb{R}\rightarrow]-1,1[$, inverse of each other, defined by%
\begin{align}
f(t)  &  =\frac{t}{\sqrt{1-t^{2}}}\qquad(-1<t<1),\label{SD06}\\
g(t)  &  =\frac{t}{\sqrt{1+t^{2}}}\qquad(t\in\mathbb{R}). \label{SD07}%
\end{align}
By the functional calculus described in Subsection \ref{FC}, we associate to
these maps the real analytic diffeomorphisms, also inverse of each other%
\begin{align*}
F  &  =\widehat{f}:\Omega\rightarrow V,\\
G  &  =\widehat{g}:V\rightarrow\Omega,
\end{align*}
where $\Omega$ is the bounded symmetric domain associated to the Jordan triple
$V$. If $\mathbf{e}=\left(  e_{1},\ldots,e_{r}\right)  $ is a frame and
$z=\sum_{j=1}^{r}\lambda_{j}e_{j}$, then%
\begin{align}
F(z)  &  =\sum_{j=1}^{r}\frac{\lambda_{j}}{\sqrt{1-\lambda_{j}^{2}}}%
e_{j}\qquad(z\in\Omega),\label{SD08}\\
G(z)  &  =\sum_{j=1}^{r}\frac{\lambda_{j}}{\sqrt{1+\lambda_{j}^{2}}}%
e_{j}\qquad(z\in V). \label{SD09}%
\end{align}
Using (\ref{P9})-(\ref{P12+}), the maps $F$ and $G$ may also be defined by%
\begin{align}
F(z)  &  =B(z,z)^{-1/4}z=\left(  \operatorname{id}_{V}-\textstyle\frac{1}{2}%
D(z,z)\right)  ^{-1/2}z\qquad(z\in\Omega),\label{SD10}\\
G(z)  &  =B(z,-z)^{-1/4}z=\left(  \operatorname{id}_{V}-\textstyle\frac{1}%
{2}D(z,-z)\right)  ^{-1/2}z\qquad(z\in V). \label{SD11}%
\end{align}

The following theorem is the main result of \cite{DiScalaLoi2006}. We give
here a different and simpler proof, using the expression of the symplectic
forms $\omega_{0}$, $\omega_{-}$, $\omega_{+}$ in generalized polar coordinates.

\begin{theorem}
\label{Duality}(Symplectic duality)
\begin{align}
F^{\ast}\omega_{0}  &  =\omega_{-},\quad F^{\ast}\omega_{+}=\omega
_{0},\label{SD03}\\
G^{\ast}\omega_{0}  &  =\omega_{+},\quad G^{\ast}\omega_{-}=\omega_{0}.
\label{SD04}%
\end{align}
\end{theorem}

\begin{proof}
In polar coordinates, the map $F$ is written%
\[
\left(  \left(  \lambda_{1},\ldots,\lambda_{r}\right)  ,\mathbf{e}\right)
\mapsto\left(  \left(  \mu_{1},\ldots,\mu_{r}\right)  ,\mathbf{e}\right)
\]
with
\begin{equation}
\mu_{j}=\frac{\lambda_{j}}{\sqrt{1-\lambda_{j}^{2}}}. \label{SD05}%
\end{equation}
As, by (\ref{S11}),%
\[
\omega_{0}=\sum_{j=1}^{r}\lambda_{j}^{2}\theta_{j0}+\sum_{\substack{j,k\\1\leq
j<k\leq r}}\left(  \lambda_{j}^{2}-\lambda_{k}^{2}\right)  \theta_{jk}%
+ \sum_{j=1}^{r}\lambda_{j}\operatorname{d}\lambda_{j}\wedge\eta_{jj},
\]
we obtain, using (\ref{SD05}),%
\begin{align*}
F^{\ast}\omega_{0}  &  =\sum_{j=1}^{r}\frac{\lambda_{j}^{2}}{1-\lambda_{j}%
^{2}}\theta_{j0}+\sum_{\substack{j,k\\1\leq j<k\leq r}}\left(  \frac
{\lambda_{j}^{2}}{1-\lambda_{j}^{2}}-\frac{\lambda_{k}^{2}}{1-\lambda_{k}^{2}%
}\right)  \theta_{jk}\\
&\quad  + \sum_{j=1}^{r}\frac{\lambda_{j}\operatorname{d}\lambda_{j}}{\left(
1-\lambda_{j}^{2}\right)  ^{2}}\wedge\eta_{jj},
\end{align*}
which, compared to (\ref{S14}), gives $F^{\ast}\omega_{0}=\omega_{-}$ on the
open dense subset $\Omega_{\mathrm{reg}}$ of regular elements, and by
continuity on all of $\Omega$.

The relation $G^{\ast}\omega_{0}=\omega_{+}$ is proved along the same lines.
The relations $F^{\ast}\omega_{+}=\omega_{0}$ and $G^{\ast}\omega_{-}%
=\omega_{0}$ follow, as $F$ and $G$ are inverse of each other. 
\end{proof}

In view of this theorem, the map $F$ (or the map $G=F^{-1}$) is called the
\emph{duality map}.

\begin{example}
(Type $I_{1,1}$) Here $V=\mathbb{C}$, $\Omega$ is the unit disc,
\[
\omega_{0}=\frac{\operatorname{i}}{2}\operatorname{d}z\wedge\operatorname{d}%
\overline{z},\quad\omega_{-}=\frac{\operatorname{i}}{2}\frac{\operatorname{d}%
z\wedge\operatorname{d}\overline{z}}{\left(  1-z\overline{z}\right)  ^{2}%
},\quad\omega_{+}=\frac{\operatorname{i}}{2}\frac{\operatorname{d}%
z\wedge\operatorname{d}\overline{z}}{\left(  1+z\overline{z}\right)  ^{2}}.
\]
The duality map is
\[
F(z)=\frac{z}{\sqrt{1-z\overline{z}}}.
\]

\end{example}

\begin{example}
(Type $I_{1,n}$) Here $V=\mathbb{C}^{n}$ with the Hermitian norm $\left\Vert
z\right\Vert ^{2}=\sum z_{j}\overline{z}_{j}$, $\Omega$ is the unit Hermitian
ball,
\[
\omega_{0}=\frac{\operatorname{i}}{2}\sum\operatorname{d}z_{j}\wedge
\operatorname{d}\overline{z}_{j},\quad\omega_{-}=\frac{\omega_{0}}{\left(
1-\left\Vert z\right\Vert ^{2}\right)  ^{n+1}},\quad\omega_{+}=\frac
{\omega_{0}}{\left(  1+\left\Vert z\right\Vert ^{2}\right)  ^{n+1}}.
\]
The duality map is%
\[
F(z)=\frac{z}{\sqrt{1-\left\Vert z\right\Vert ^{2}}}.
\]
\end{example}

\section{The bisymplectomorphism group\label{SEC4}}

\subsection{}

In this section we study to which extent a diffeomorphism $F:\Omega\rightarrow
V$ satisfying the property (\ref{SD03}) is unique. If two diffeomorphisms
$F_{1},F_{2}:\Omega\rightarrow V$ satisfy (\ref{SD03}), then $F_{2}=F_{1}\circ
f$, where $f:\Omega\rightarrow\Omega$ preserves $\omega_{0}$ and $\omega_{-}$.
This leads us to the following definition.

\begin{definition}
A \emph{bisymplectomorphism of} $\Omega$ is a diffeomorphism $f:\Omega
\rightarrow\Omega$ which satisfies
\begin{align}
f^{\ast}\omega_{0}  &  =\omega_{0},\label{G06}\\
f^{\ast}\omega_{-}  &  =\omega_{-}. \label{G07}%
\end{align}

\end{definition}

Clearly, bisymplectomorphisms of $\Omega$ form a group, which will be denoted
by $\mathcal{B}(\Omega)$ and called the \emph{bisymplectomorphism group} of
$\Omega$.

From (\ref{G05}), we derive a characterization of $\mathcal{B}(\Omega)$ in
terms of the Bergman operator:

\begin{proposition}
\label{general} Let $\Omega$ be a bounded symmetric domain. Then a
diffeomorphism $f\in\operatorname{Diff}(\Omega)$ is a bisymplectomorphism if
and only if it satisfies%
\begin{align}
f^{\ast}\omega_{0}  &  =\omega_{0},\label{G08}\\
B\left(  f(z),f(z)\right)  \circ\operatorname{d}f(z)  &  =\operatorname{d}%
f(z)\circ B(z,z)\quad(z\in\Omega). \label{G09}%
\end{align}

\end{proposition}

Note that the second condition implies that the tangent map $\operatorname{d}%
f(z)$ maps invariant subspaces of $B_{z}=B(z,z)$ to invariant subspaces of
$B_{f(z)}$, and that $B_{f(z)}$ has the same eigenvalues as $B_{z}$.

\begin{proof}
Let $z\in\Omega$, $u,v\in T_{z}\Omega$. We have from (\ref{G05}) and
(\ref{G06})%
\begin{align*}
\omega_{-}(z)(u,v)  &  =\omega_{0}(B(z,z)^{-1}u,v)\\
&  =\omega_{0}(\operatorname{d}f(z)B(z,z)^{-1}u,\operatorname{d}f(z)v),\\
\left(  f^{\ast}\omega_{-}\right)  (z)(u,v)  &  =\omega_{-}%
(f(z))(\operatorname{d}f(z)u,\operatorname{d}f(z)v)\\
&  =\omega_{0}\left(  B(f(z),f(z))^{-1}\operatorname{d}f(z)u,\operatorname{d}%
f(z)v\right)  ,
\end{align*}
so that, assuming (\ref{G06}), the condition (\ref{G07}) is equivalent to
\[
\omega_{0}(\operatorname{d}f(z)B(z,z)^{-1}u,\operatorname{d}f(z)v)=\omega
_{0}\left(  B(f(z),f(z))^{-1}\operatorname{d}f(z)u,\operatorname{d}%
f(z)v\right)
\]
for all $u,v$. As $\omega_{0}$ is non singular and $\operatorname{d}f(z)$ is
bijective, this is equivalent to%
\[
\operatorname{d}f(z)\circ B(z,z)^{-1}=B(f(z),f(z))^{-1}\circ\operatorname{d}%
f(z),
\]
that is, to (\ref{G09}). 
\end{proof}

\subsection{}

Here we study diffeomorphisms of $\Omega$ satisfying the condition
(\ref{G09}). Recall that $B_{z}=B(z,z):V\rightarrow V$ is a $\mathbb{C}%
$-linear operator, self-adjoint w.r.~to the Hermitian metric $h_{0}$, positive
if $z\in\Omega$. Let $r$ denote the \emph{rank} of $\Omega$ and $V$.

For $z\in\Omega$, consider the spectral decomposition
\begin{equation}
z=\lambda_{1}e_{1}+\lambda_{2}e_{2}+\cdots+\lambda_{s}e_{s},\quad1>\lambda
_{1}>\lambda_{2}>\cdots>\lambda_{s}>0, \label{G17}%
\end{equation}
where $s=\operatorname{rk}z\leq r=\operatorname{rk}V$. An element is called
\emph{regular} if $\operatorname{rk}z=\operatorname{rk}V$; for regular
elements, which form an open dense subset of $\Omega$, the decomposition
(\ref{G17}) is the decomposition using generalized polar coordinates.

Let
\[
V={\displaystyle\bigoplus\limits_{0\leq i\leq j\leq s}} V_{ij}%
\]
be the simultaneous Peirce decomposition relative to $\left(  e_{1}%
,\ldots,e_{s}\right)  $. Note that some subspaces $V_{ij}$ may be $0$. The
operator $B(z,z)$ may only have the eigenvalues
\begin{align}
(1-\lambda_{i}^{2})^{2}\quad(1  &  \leq i\leq s),\quad(1-\lambda_{i}%
^{2})(1-\lambda_{j}^{2})\quad(1\leq i<j\leq s),\label{G10}\\
(1-\lambda_{i}^{2})\quad(1  &  \leq i\leq s),\quad1,\nonumber
\end{align}
which occur respectively on the subspaces
\begin{equation}
V_{ii},\quad V_{ij},\quad V_{0i},\quad V_{00}. \label{G11}%
\end{equation}

The relation (\ref{G09}) then implies that $B(z,z)$ and $B(f(z),f(z))$ have
the same eigenvalues with the same multiplicities. Moreover,\emph{ if all
eigenvalues in the list (\ref{G10}) occurring for non-zero }$V_{ij}$\emph{ are
different}, the non-zero subspaces of the list (\ref{G11}) are the eigenspaces
of $B(z,z)$ and are mapped by $\operatorname{d}f(z)$ to the corresponding
eigenspaces of $B(f(z),f(z))$.

\begin{definition}
\label{super} Let $\Omega$ be an irreducible bounded symmetric domain and
denote by $\left(  V,\left\{  ~,~,~\right\}  \right)  $ be the corresponding
Jordan triple. An element $z\in V$ is called \emph{super-regular} if $z$ is
regular and if all eigenvalues in the list (\ref{G10}), occurring on non-zero
$V_{ij}$, are different.
\end{definition}

If $V$ is of tube type or of rank $1$, any regular element $z\in V$ is
super-regular. If $V$ is not of tube type, an element $z\in V$ is
super-regular if it is regular and if its spectral values satisfy
\begin{equation}
1-\lambda_{i}^{2}\neq\left(  1-\lambda_{j}^{2}\right)  ^{2} \label{G12}%
\end{equation}
for all $(i,j)$, $i<j$. Clearly, super-regular elements form an open dense
subset of $V$.

Let $z\in\Omega$ be an element of \emph{rank one}. Then the spectral
decomposition of $z$ is
\[
z=\lambda_{1}e_{1},\quad0<\lambda_{1}<1,
\]
where $e_{1}$ is a tripotent; the associated Peirce decomposition is
\[
V_{11}=V_{2}(e_{1}),\quad V_{01}=V_{1}(e_{1}),\quad V_{00}=V_{0}(e_{1})
\]
and the eigenvalues of $B(z,z)$ on these subspaces are respectively%
\[
\left(  1-\lambda_{1}^{2}\right)  ^{2},\quad1-\lambda_{1}^{2},\quad1.
\]
It then follows that
\[
f(z)=\lambda_{1}\varepsilon_{1},
\]
where $\varepsilon_{1}$ is a tripotent such that $V_{2}(\varepsilon
_{1})=\operatorname{d}f(z)V_{2}(e_{1})$, which means that $e_{1}$ and
$\varepsilon_{1}$ have the same height $\dim V_{2}(e_{1})=\dim V_{2}%
(\varepsilon_{1})$. In particular, if $e_{1}$ is minimal (resp. maximal), then
$\varepsilon_{1}$ is minimal (resp. maximal).

Let $V_{z}$ be the $\mathbb{C}$-subspace generated by the odd powers
$z,\ldots,z^{(2p+1)},\ldots$ Then $\dim_{\mathbb{C}}V_{z}=\operatorname{rk}%
z\leq r$, and $\dim_{\mathbb{C}}V_{z}=r$ if and only if $z$ is a regular
element. Denote
\[
P_{z}=V_{z}\cap\Omega.
\]
For $z\in\Omega$, $z\neq0$ with spectral decomposition $z=\sum_{j=1}^{k}%
\alpha_{j}e_{j}$ ($\alpha_{1}>\cdots>\alpha_{k}>0$; $k=\operatorname{rk}_{V}%
x$), $P_{z}$ is the $k$-polydisc%
\[
P_{z}=\left\{  u=\sum u_{j}e_{j}\mid\left\vert u_{j}\right\vert <1\right\}  .
\]
Note that if $u\in P_{z}$ and $\operatorname{rk}_{V}u=\operatorname{rk}_{V}z$,
then $V_{u}=V_{z}$ and $P_{u}=P_{z}$.

\begin{lemma}
Let $f:\Omega\rightarrow\Omega$ be a diffeomorphism of $\Omega$ such that
\[
B\left(  f(z),f(z)\right)  \circ\operatorname{d}f(z)=\operatorname{d}f(z)\circ
B(z,z)
\]
for all $z\in\Omega$. Then for all $z\in\Omega$,
\begin{equation}
\operatorname{d}f(z)V_{z}=V_{f(z)}. \label{G15}%
\end{equation}

\end{lemma}

\begin{proof}
Let $z$ be a super-regular element and let%
\[
z=\alpha_{1}e_{1}+\cdots+\alpha_{r}e_{r}%
\]
be the spectral decomposition of $z$ in $V$; then $V_{z}:=\mathbb{C}%
e_{1}\oplus\cdots\oplus\mathbb{C}e_{r}$ is the sum of the eigenspaces of
$B(z,z)$ relative to the eigenvalues%
\[
\left(  1-\alpha_{j}^{2}\right)  ^{2}\quad(1\leq j\leq r),
\]
so that it follows from (\ref{G09}) that%
\[
\operatorname{d}f(z)V_{z}=V_{f(z)}.
\]
By continuity, (\ref{G15}) holds for all $z\in\Omega$. 
\end{proof}

\begin{proposition}
\label{radiality}Let $f:\Omega\rightarrow\Omega$ be a diffeomorphism of
$\Omega$ such that
\[
B\left(  f(z),f(z)\right)  \circ\operatorname{d}f(z)=\operatorname{d}f(z)\circ
B(z,z)
\]
for all $z\in\Omega$. Then for any element $z\in\Omega$, $z\neq0$, we have%
\[
f\left(  P_{z}\right)  =P_{f(z)}.
\]
\end{proposition}

\begin{proof}
Let $s=\operatorname{rk}_{V}z$. We already know that $\operatorname{d}%
f(z)V_{z}=V_{f(z)}$, which implies that
$\operatorname{rk}_{V}f(z)=\operatorname{rk}_{V}z$. So
there exist continuous functions $\beta_{j}$ such that
\[
\operatorname{d}f(z)u=\sum_{j=1}^{s}\beta_{j}(u,z)\left(  f(z)\right)
^{(2j+1)}\quad(z\in\Omega,\ u\in P_{z}).
\]
Consider a $C^{1}$ path $\eta:[0,1]\rightarrow P_{z}$ from $\eta(0)=z$ to
$\eta(1)=w$; let $g(t)=f(\eta(t))$. Then $g$ satisfies the differential
equation%
\begin{align}
g^{\prime}(t)  &  =\sum_{j=1}^{s}\beta_{j}(\eta^{\prime}(t),\eta(t))\left(
g(t)\right)  ^{(2j+1)},\label{G16}\\
g(0)  &  =f(z).\nonumber
\end{align}
Let
\[
f(z)=\sum_{k=1}^{s}\alpha_{k}\varepsilon_{k}%
\]
be the spectral decomposition of $f(z)$. Let $h:[0,1]\rightarrow\mathbb{C}%
^{s}$ be the solution of the differential system%
\begin{align*}
h_{k}^{\prime}(t)  &  =\sum_{j=1}^{s}\beta_{j}(\eta^{\prime}(t),\eta
(t))\left(  h_{k}(t)\right)  ^{2j+1)}\quad(1\leq k\leq s),\\
h_{k}(0)  &  =\alpha_{k}.
\end{align*}
Then the solution of (\ref{G16}) is
\[
g(t)=\sum_{k=1}^{s}h_{k}(t)\varepsilon_{k},
\]
which shows that $g(1)=f(w)$ belongs to $P_{f(z)}$. 
\end{proof}

Let $K$ denote the group of linear automorphisms of $\Omega$.

\begin{proposition}
Let $f:\Omega\rightarrow\Omega$ be a diffeomorphism of $\Omega$ such that
\[
B\left(  f(z),f(z)\right)  \circ\operatorname{d}f(z)=\operatorname{d}f(z)\circ
B(z,z)
\]
for all $z\in\Omega$. Then any $K$-orbit is globally invariant by $f$.
\end{proposition}

\begin{proof}
Assume first that $z=\lambda_{1}e_{1}+\cdots+\lambda_{r}e_{r}$ is super
regular. Then $f(z)=\mu_{1}\varepsilon_{1}+\cdots+\mu_{r}\varepsilon_{r}$, and
equality between the eigenvalues of $B(z,z)$ and $B(f(z),f(z))$ implies
$\lambda_{j}=\mu_{j}$, hence $f(z)\in Kz$ and $f(Kz)=Kz$. By continuity, this
also holds for any $z\in\Omega$. 
\end{proof}

\subsection{}

We now go back to bisysmplectomorphisms of $\Omega$.

\begin{proposition}
\label{diff_at_0}Let $\Omega$ be a bounded circled symmetric domain and denote
by $K$ the group of linear automorphisms of $\Omega$. For each
bisymplectomorphism $f\in\mathcal{B}(\Omega)$, we have $\operatorname{d}%
f(0)\in K$.
\end{proposition}

As $K\subset\mathcal{B}(\Omega)$, it will be sufficient to study the subgroup
\begin{equation}
\mathcal{B}_{0}(\Omega)=\left\{  f\in\mathcal{B}(\Omega)\mid\operatorname{d}%
f(0)=\operatorname{id}_{V}\right\}  . \label{G14}%
\end{equation}

\begin{proof}
Let $\left(  e_{1},\ldots,e_{r}\right)  $ be a frame of $V$ and let
$V={\displaystyle\bigoplus}V_{ij}$ be the associated simultaneous Peirce
decomposition. Consider a regular element
\[
z=\alpha_{1}e_{1}+\cdots+\alpha_{r}e_{r},
\]
$1>\alpha_{1}>\cdots>\alpha_{r}>0$. For $f\in\mathcal{B}(\Omega)$ and $z$
super-regular, we have $f(z)=\alpha_{1}\varepsilon_{1}+\cdots+\alpha
_{r}\varepsilon_{r}$, where $\left(  \varepsilon_{1},\ldots,\varepsilon
_{r}\right)  $ is a frame of $V$, which may depend on $\left(  \alpha
_{1},\ldots,\alpha_{r}\right)  $. As $V_{jj}=\mathbb{C}e_{j}$ is the
eigenspace of $B(z,z)$ for the eigenvalue $\left(  1-\alpha_{j}^{2}\right)
^{2}$, we deduce from (\ref{G09}) that
\[
\operatorname{d}f(z)\left(  e_{j}\right)  \in\mathbb{C\varepsilon}_{j}.
\]
By continuity, there exists a frame $\left(  e_{1}^{\prime},\ldots
,e_{r}^{\prime}\right)  $ such that
\[
\operatorname{d}f(0)\left(  e_{j}\right)  \in\mathbb{C}e_{j}^{\prime}%
\qquad(1\leq j\leq r).
\]
Multiplying $e_{j}^{\prime}$ by a suitable complex of modulus $1$, we may
assume that
\begin{equation}
\operatorname{d}f(0)\left(  e_{j}\right)  =\lambda_{j}^{\prime}e_{j}^{\prime
},\quad\lambda_{j}^{\prime}>0\quad(1\leq j\leq r). \label{G13}%
\end{equation}
We have also
\[
\operatorname{d}f(0)(e_{j})=\lim_{t\rightarrow0+}\frac{f(te_{j})}{t}%
\]
and $f(te_{j})$ is a multiple of a minimal tripotent, with spectral norm
$\left\Vert f(te_{j})\right\Vert \leq t\left\Vert e_{j}\right\Vert $, hence
$\left\Vert \operatorname{d}f(0)(e_{j})\right\Vert \leq\left\Vert
e_{j}\right\Vert $ and $\lambda_{j}^{\prime}\leq1$. The same argument applied
to $f^{-1}$ gives $\lambda_{j}^{\prime}=1$. So the image of a frame $\left(
e_{1},\ldots,e_{r}\right)  $ under $\operatorname{d}f(0)$ is a frame $\left(
e_{1}^{\prime},\ldots,e_{r}^{\prime}\right)  $ and the Peirce spaces $V_{jk}$
relative to $\left(  e_{1},\ldots,e_{r}\right)  $ are mapped by
$\operatorname{d}f(0)$ onto the corresponding Peirce spaces $V_{jk}^{\prime}$
relative to $\left(  e_{1}^{\prime},\ldots,e_{r}^{\prime}\right)  $. In
particular, $\operatorname{d}f(0)$ is an $\mathbb{R}$-linear map from
$\mathbb{C}e_{j}$ onto $\mathbb{C}e_{j}^{\prime}$. We have $\operatorname{d}%
f(0)(\operatorname{i}e_{j})=\beta e_{j}^{\prime}$, with $\left\vert
\beta\right\vert =1$. As
\[
\omega_{0}(u,v)=\textstyle\frac{\operatorname{i}}{2}\left(  \operatorname{tr}%
D(u,v)-\operatorname{tr}D(v,u)\right)  ,
\]
we have
\begin{align*}
\omega_{0}(e_{j},\operatorname{i}e_{j})  &  =\operatorname{tr}D(e_{j}%
,e_{j}),\\
\omega_{0}\left(  \operatorname{d}f(0)\left(  e_{j}\right)  ,\operatorname{d}%
f(0)\left(  \operatorname{i}e_{j}\right)  \right)   &  =\operatorname{Im}%
\beta\operatorname{tr}D(e_{j}^{\prime},e_{j}^{\prime}).
\end{align*}
From $\operatorname{tr}D(e_{j}^{\prime},e_{j}^{\prime})=\operatorname{tr}%
D(e_{j},e_{j})$ and $f^{\ast}\omega_{0}=\omega_{0}$, we get $\operatorname{Im}%
\beta=1$, $\beta=\operatorname{i}$, which means that $\operatorname{d}f(0)$ is
$\mathbb{C}$-linear on $\mathbb{C}e_{1}\oplus\cdots\oplus\mathbb{C}e_{r}$.
Finally, $\operatorname{d}f(0)$ is $\mathbb{C}$-linear on $V$, maps minimal
tripotents to minimal tripotents and $\Omega$ to $\Omega$. 
\end{proof}

\subsection{The unit disc}

Let $\Delta$ be the unit disc of $\mathbb{C}$. The associated triple product
is $\left\{  u,v,w\right\}  =2u\overline{v}w$. The Bergman operator is
$B(z,z)w=(1-|z|^{2})^{2}w$. The Bergman metric is
\[
h_{z}(u,v)=\frac{u\overline{v}}{(1-|z|^{2})^{2}}.
\]
The two symplectic forms are
\begin{align}
\omega_{0}  &  ={\textstyle\frac{\operatorname{i}}{2}}\operatorname{d}z\wedge\operatorname{d}%
\overline{z},\label{disc5}\\
\omega_{-}  &  =\frac{\omega_{0}}{(1-|z|^{2})^{2}}. \label{disc6}%
\end{align}
Denote by $S^{1}$ the unit circle in $\mathbb{C}$ and consider the
\textquotedblleft polar coordinates\textquotedblright\ diffeomorphism
\begin{align*}
\Theta:(0,1)\times S^{1}  &  \rightarrow\Delta\setminus\{0\},\\
\left(  r,\zeta\right)   &  \mapsto r\zeta.
\end{align*}
Then we have
\[
\Theta^{\ast}\omega_{0}=r\operatorname{d}r\wedge\frac{\operatorname{d}\zeta
}{\operatorname{i}\zeta}.
\]

The following theorem characterizes the elements of $\mathcal{B}(\Delta)$.

\begin{theorem}
\label{rk1-1}The elements $f\in\mathcal{B}(\Delta)$ are the maps defined by
\[
f(z)=u\left(  |z|^{2}\right)  \,z\qquad\left(  z\in\Delta\right)  ,
\]
where $u$ is a smooth function $u:[0,1)\rightarrow S^{1}\simeq U(1)$.
\end{theorem}

In other words, the restriction of $f$ to a circle of radius $r$ ($0<r<1$) is
the rotation $u\left(  r^{2}\right)  $.

\begin{proof}
From (\ref{disc6}), we see that a diffeomorphism $f:\Delta\rightarrow\Delta$
is a bisymplectomorphism if and only if $f$ preserves $\omega_{0}$ and
$|f(z)|=|z|$ for all $z\in\Delta$.

If $|f(z)|=|z|$ for all $z\in\Delta$, the map $F=\Theta^{-1}\circ f\circ
\Theta$ may be written
\[
(r,\zeta)\overset{F}{\rightarrow}(r,Z(r,\zeta))
\]
for some smooth function $Z:(0,1)\times S^{1}\rightarrow S^{1}$. We have
\[
F^{\ast}\left(  \Theta^{\ast}\omega_{0}\right)  = r\operatorname{d}%
r\wedge\frac{\operatorname{d}Z}{\operatorname{i}Z}= r\operatorname{d}%
r\wedge\frac{\operatorname{d}_{\zeta}Z}{\operatorname{i}Z}.
\]
If $f$ preserves $\omega_{0}$, then $F$ preserves $\Theta^{\ast}\omega_{0}$,
which implies%
\begin{equation}
\frac{\operatorname{d}_{\zeta}Z}{Z}=\frac{\operatorname{d}\zeta}{\zeta}.
\label{disc1}%
\end{equation}
For $r$ fixed, let $u_{r}(\zeta)=Z(r,\zeta)$ ; the condition (\ref{disc1}) is
then equivalent to
\begin{equation}
\frac{\operatorname{d}u_{r}}{u_{r}}=\frac{\operatorname{d}\zeta}{\zeta}
\label{disc2}%
\end{equation}
for all $r\in(0,1)$. The condition (\ref{disc2}) is in turn equivalent to
\begin{equation}
u_{r}(\zeta)=v(r)\zeta, \label{disc3}%
\end{equation}
with $v(r)\in S^{1}$. The function $v$, which is given by%
\[
v(r)=\frac{Z(r,\zeta)}{\zeta},
\]
is smooth on $(0,1)$ and
\[
f(r\zeta)=rv(r)\zeta,
\]
or
\begin{equation}
f(z)=v\left(  \left\vert z\right\vert \right)  z \label{disc4}%
\end{equation}
for $z\in\Delta$, $z\neq0$; on the other hand, $f(0)=0$. Let $g$ be the
restriction of $f$ to $(-1,1)$; from (\ref{disc4}), we see that $g$ is odd.
For $x\neq0,$ we have
\[
g(x)=\int_{0}^{1}\frac{\operatorname{d}}{\operatorname{d}s}%
g(sx)\operatorname{d}s=x\int_{0}^{1}g^{\prime}(sx)\operatorname{d}s,
\]
which shows that
\[
v(r)=\int_{0}^{1}g^{\prime}(sr)\operatorname{d}s
\]
extends to a smooth even function $v:(-1,1)\rightarrow S^{1}$. Let
$u:[0,1)\rightarrow S^{1}$ be defined by
\[
u(r)=v\left(  \sqrt{r}\right)  .
\]
It follows then from Whitney's theorem \cite{Whitney1943} that $u$ is smooth
at $0$.

Conversely, if
\[
f(z)=u(\left\vert z\right\vert ^{2})z,
\]
with $u:[0,1)\rightarrow S^{1}$ smooth, then $f$ satisfies $\left\vert
f(z)\right\vert =\left\vert z\right\vert $, $f$ is a diffeomorphism with
inverse $f^{-1}(w)=u(\left\vert w\right\vert ^{2})^{-1}w$; also, $f$ preserves
$\omega_{0}$ on $\Delta\setminus\{0\}$, hence on $\Delta$ by continuity. This
implies that $f\in\mathcal{B}(\Delta)$. 
\end{proof}

\subsection{The polydisc}

Let $\Delta^{r}\subset\mathbb{C}^{r}$ be the product of $r$ unit discs. The
Jordan triple product on $V=\mathbb{C}^{r}$ is just the component-wise product%
\[
\left\{  x,y,z\right\}  =2\left(  x_{1}\overline{y_{1}}z_{1},\ldots
,x_{r}\overline{y_{r}}z_{r}\right)  .
\]

Let $\left(  e_{1},\ldots,e_{r}\right)  $ denote the canonical basis of
$\mathbb{C}^{r}$. The minimal tripotents of $V$ are the elements $\lambda
_{j}e_{j}$, $1\leq j\leq r$, $\left\vert \lambda_{j}\right\vert =1$. Any frame
(maximal ordered set of mutually orthogonal tripotents) has the form%
\[
\left(  \lambda_{j}e_{\sigma(j)}\right)  _{1\leq j\leq r},
\]
where $\left\vert \lambda_{j}\right\vert =1$ and $\sigma\in\mathfrak{S}_{r}$
is a permutation of $\{1,\cdots,r\}$. The corresponding Peirce decomposition
is $V=\mathbb{C}e_{\sigma(1)}\oplus\cdots\oplus\mathbb{C}e_{\sigma(r)}$.

\begin{theorem}
\label{polydisc}A diffeomorphism $f:\Delta^{r}\rightarrow\Delta^{r}$ belongs
to $\mathcal{B}_{0}\left(  \Delta^{r}\right)  $ if and only if there exist
smooth functions $u_{j}:[0,1)\rightarrow S^{1}$ such that $u_{j}(0)=1$ and
\begin{equation}
f\left(  z_{1},\ldots,z_{r}\right)  =\sum_{j=1}^{r}u_{j}\left(  |z_{j}%
|^{2}\right)  \,z_{j}e_{j}\qquad\left(  z_{j}\in\Delta\right)  .
\label{polydisc1}%
\end{equation}
\end{theorem}

\begin{proof}
Let $f\in\mathcal{B}_{0}(\Delta^{r})$ be a bisymplectomorphism with
$\operatorname{d}f(0)=\operatorname{id}_{V}$. Consider a regular element
$z\in\Delta^{r}$, that is,
\[
z=z_{1}e_{1}+\cdots+z_{r}e_{r},
\]
with all $\left\vert z_{j}\right\vert $ different. The spaces $\mathbb{C}%
e_{j}$ of the corresponding Peirce decomposition are mapped by
$\operatorname{d}f(z)$ to the spaces $\mathbb{C}e_{\sigma(j)}$ of another
Peirce decomposition, for some permutation $\sigma\in\mathfrak{S}_{r}$. This
means that $\left[  \operatorname{d}f(z)\left(  e_{j}\right)  \right]
=\left[  e_{\sigma(j)}\right]  $ for all regular $z\in\Delta^{r}$, where
$\left[  \quad\right]  $ denotes the class in $\mathbb{P}(\mathbb{C}^{r})$; by
continuity, this is true for all $z\in\Delta^{r}$. As $\left\{  \left[
e_{1}\right]  ,\ldots,\left[  e_{r}\right]  \right\}  $ is discrete in
$\mathbb{P}(\mathbb{C}^{r})$ and $\operatorname{d}f(0)=\operatorname{id}$, we
have $\left[  \operatorname{d}f(z)\left(  e_{j}\right)  \right]  =\left[
e_{j}\right]  $ and
\begin{equation}
\operatorname{d}f(z)\left(  \mathbb{C}e_{j}\right)  =\mathbb{C}e_{j}
\label{polydisc2}%
\end{equation}
for all $z\in\Delta^{r}$. This shows that $f\left(  z_{1},\ldots,z_{r}\right)
=\sum_{j=1}^{r}f_{j}(z)e_{j}$, with $\left\vert f_{j}(z)\right\vert
=\left\vert z_{j}\right\vert $. From (\ref{polydisc2}), we deduce that $f_{j}$
depends only of $z_{j}$. Each $f_{j}$ is then a bisymplectomorphism of the
unit disc. According to Theorem \ref{rk1-1}, there exists a smooth function
$u_{j}:[0,1)\rightarrow S^{1}\simeq U(1)$ such that%
\[
f_{j}(z_{j})=u_{j}\left(  |z_{j}|^{2}\right)  \,z_{j}\qquad\left(  z_{j}%
\in\Delta\right)  .
\]
Finally, any $f\in\mathcal{B}_{0}(\Delta^{r})$ has the form (\ref{polydisc1}).
Conversely, each $f$ of the form (\ref{polydisc1}) is easily seen
to be a bisymplectomorphism.
 \end{proof}

\subsection{The general case}

We assume now that the domain $\Omega$ is \emph{irreducible}, that is, not a
product of two bounded symmetric domains. In Theorem \ref{structure}, we will
characterize the bisymplectomorphisms of $\Omega$. This theorem may be
considered as a kind of Schwarz lemma.

Let $f\in\mathcal{B}_{0}(\Omega)$. If $z=\sum\lambda_{j}e_{j}$ is a regular
element, then $P_{z}$ is the polydisc $(  \mathbb{C}e_{1}\oplus\cdots\oplus\mathbb{C}%
e_{r})  \cap\Omega$ and $f\left(  P_{z}\right)  =P_{z}$. For a frame
$\mathbf{e}=\left(  e_{1},\ldots,e_{r}\right)  $, let $P(\mathbf{e)=}\left(
\mathbb{C}e_{1}\oplus\cdots\oplus\mathbb{C}e_{r}\right)  \cap\Omega$. By the
same argument as in the proof of Theorem \ref{polydisc}, the restriction of
$f$ to $P(\mathbf{e)}$ has the form
\[
\sum_{j=1}^{r}\lambda_{j}e_{j}\mapsto\sum\lambda_{j}u_{j}\left(  |\lambda
_{j}|^{2}\right)  e_{j},
\]
where the $u_{j}$'s are smooth functions $u_{j}:[0,1)\rightarrow S^{1}\simeq
U(1)$. In the \textquotedblleft polar coordinates\textquotedblright\ $\left(
\left(  \lambda_{1},\ldots,\lambda_{r}\right)  ,\mathbf{e}=\left(
e_{1},\ldots,e_{r}\right)  \right)  $, the map $f$ is then represented by
\begin{equation}
\left(  \left(  \lambda_{1},\ldots,\lambda_{r}\right)  ,\mathbf{e}\right)
\mapsto\left(  \left(  \lambda_{1},\ldots,\lambda_{r}\right)  ,\left(
u_{1}\left(  \lambda_{1}^{2},\mathbf{e}\right)  e_{1},\ldots,u_{r}\left(
\lambda_{r}^{2},\mathbf{e}\right)  e_{r}\right)  \right)  . \label{G19}%
\end{equation}
Let $w_{j}\left(  \lambda_{j},\mathbf{e}\right)  =u_{j}\left(  \lambda_{j}%
^{2},\mathbf{e}\right)  $. We obtain%
\begin{align*}
f^{\ast}\eta_{jj}  &  ={\textstyle\frac{\operatorname{i}}{2}}f^{\ast}\sum_{m=1}^{n}e_{jm}%
\operatorname{d}\overline{e}_{jm}={\textstyle\frac{\operatorname{i}}{2}}\sum_{m=1}^{n}w_{j}%
e_{jm}\left(  \overline{w}_{j}\operatorname{d}\overline{e}_{jm}+\overline
{e}_{jm}\operatorname{d}\overline{w}_{j}\right) \\
&  =\eta_{jj}-{\textstyle\frac{\operatorname{i}}{2}}\left(  e_{j}\mid e_{j}\right)  \frac
{\operatorname{d}w_{j}}{w_{j}},\\
f^{\ast}\omega_{jj}  &  =\omega_{jj}.
\end{align*}
As $\Omega$ is irreducible, $\left(  e_{j}\mid e_{j}\right)  $ has the same
value $g$ for all minimal tripotents. Finally%
\begin{align*}
f^{\ast}\omega_{0}  &  =\sum_{j=1}^{r}\lambda_{j}^{2}\omega_{jj}+ \sum
_{j=1}^{r}\lambda_{j}\operatorname{d}\lambda_{j}\wedge\left(  \eta
_{jj}-\operatorname{i}g\frac{\operatorname{d}w_{j}}{w_{j}}\right) \\
&  =\omega_{0}- \operatorname{i}g\sum_{j=1}^{r}\lambda_{j}\operatorname{d}%
\lambda_{j}\wedge\frac{\operatorname{d}w_{j}}{w_{j}}.
\end{align*}
Then $f^{\ast}\omega_{0}=\omega_{0}$ implies that $f\in\mathcal{B}_{0}%
(\Omega)$, written in the form (\ref{G19}), satisfies%
\begin{equation}
\sum_{j=1}^{r}\lambda_{j}\operatorname{d}\lambda_{j}\wedge\frac
{\operatorname{d}w_{j}}{w_{j}}=0. \label{G18}%
\end{equation}
As $w_{j}$ depends only on $\lambda_{j}$ and $\mathbf{e}$, this implies that
$\operatorname{d}_{\mathbf{e}}w_{j}=0$. As the manifold of frames is connected
when the domain $\Omega$ is irreducible, $w_{j}$ does not depend on
$\mathbf{e}\in\mathcal{F}$. As a permutation of a frame is again a frame, we
have $w_{1}=\cdots=w_{r}$ and $u_{1}=\cdots=u_{r}$.

Finally, an element $f\in\mathcal{B}_{0}(\Omega)$ is written in polar
coordinates%
\begin{equation}
\left(  \left(  \lambda_{1},\ldots,\lambda_{r}\right)  ,\left(  e_{1}%
,\ldots,e_{r}\right)  \right)  \mapsto\left(  \left(  \lambda_{1}%
,\ldots,\lambda_{r}\right)  ,\left(  u\left(  \lambda_{1}^{2}\right)
e_{1},\ldots,u\left(  \lambda_{r}^{2}\right)  e_{r}\right)  \right)  ,
\label{G20}%
\end{equation}
where $u$ is a smooth function $u:[0,1)\rightarrow S^{1}\simeq U(1)$.

\begin{theorem}
\label{structure}Let $\Omega$ be an irreducible Hermitian bounded circled
symmetric domain and let $K$ be the isotropy group of $0$. The analytic (resp.
$C^{\infty}$) bisymplectomorphisms of $\Omega$ are the maps $\phi=f\circ g$,
where $g=\operatorname{d}\phi(0)\in K$ and $f$ is associated to $v(t)=tu(t^{2}%
)$, with $u:[0,1)\rightarrow S^{1}\simeq U(1)$ analytic (resp. $C^{\infty}$)
and $u(0)=1$.
\end{theorem}

\end{document}